\documentclass[11pt]{article}

\usepackage{amssymb,amsmath,latexsym}
\usepackage[dvips]{color}

\begin{document}

\newtheorem{thm}{Theorem}[section]
\newtheorem{lemma}[thm]{Lemma}
\newtheorem{defn}{Definition}[section]
\newtheorem{prop}[thm]{Proposition}
\newtheorem{corollary}[thm]{Corollary}
\newtheorem{remark}[thm]{Remark}
\newtheorem{example}[thm]{Example}
\newtheorem{proof}{proof}
\newtheorem{assumption}{Assumption}
\numberwithin{equation}{section}

\def\ee{\varepsilon}
\def\qed{{\hfill $\Box$ \bigskip}}
\def\MM{{\cal M}}
\def\BB{{\cal B}}
\def\LL{{\cal L}}
\def\FF{{\cal F}}
\def\EE{{\cal E}}
\def\QQ{{\cal Q}}
\def\AA{{\cal A}}

\def\R{{\mathbb R}}
\def\N{{\mathbb N}}
\def\E{{\bf E}}
\def\F{{\bf F}}
\def\H{{\bf H}}
\def\P{{\bf P}}
\def\Q{{\bf Q}}
\def\S{{\bf S}}
\def\J{{\bf J}}
\def\K{{\bf K}}
\def\F{{\bf F}}
\def\A{{\bf A}}
\def\loc{{\bf loc}}
\def\eps{\varepsilon}
\def\semi{{\bf semi}}
\def\wh{\widehat}
\def\pf{\noindent{\bf Proof} }
\def\dim{{\rm dim}}

\title{\Large \bf  Strong law of large number of a class of super-diffusions}
\author{ \bf Rong-Li Liu, \hspace{1mm}\hspace{1mm} Yan-Xia
Ren\footnote{The research of this author is supported by NSFC (Grant
No.  10871103 and 10971003) and Specialized Research Fund for the
Doctoral Program of Higher Education.\hspace{1mm} }
\hspace{1mm}\hspace{1mm} and Renming Song\hspace{1mm} }

 \date{}\maketitle

 {\narrower{\narrower

\centerline{\bf Abstract}

\bigskip
In this paper we prove that, under certain conditions, a strong law
of large numbers holds for a class of super-diffusions $X$
corresponding to the evolution equation $\partial_t u_t=L u_t+\beta
u_t-\psi(u_t)$ on a bounded domain $D$ in $\R^d$,
where $L$ is the generator of the underlying diffusion and the
branching mechanism
 $\psi(x,\lambda)=\frac{1}{2}\alpha(x)\lambda^2+\int_0^\infty
(e^{-\lambda r}-1+\lambda r)n(x, {\rm d}r)$ satisfies $\sup_{x\in
D}\int_0^\infty (r\wedge r^2) n(x,{\rm d}r)<\infty$.

\bigskip
\noindent{{\bf Keywords\hspace{2mm}} Super-diffusion, martingale,
point process, principal eigenvalue, strong law of large numbers}
\medskip

\noindent{{\bf  2010 MR Subject Classification}\hspace*{10pt}
 %60J80,\hspace{3pt} 60G57,\hspace{3pt} 60G17}
 60J68,\hspace{3pt} 60G55,\hspace{3pt} 60G57}

\par}\par}

\bigskip
\begin{doublespace}
\parskip -0.5mm \vspace{2 true mm}
\section{Introduction }
\subsection{Motivation}

Recently many people (see \cite{CRW, CS, E,EHK,ET,EW,W} and the
references therein) have studied limit theorems for branching Markov
processes or super-processes using the principal eigenvalue and
ground state of the linear part of the characteristic equations. All
the papers above, except \cite{ET}, assumed that
the branching mechanisms satisfy a second moment condition.
In \cite{ET}, a $(1+\theta)$-moment condition, $\theta>0$, on the
branching mechanism is assumed instead.

In \cite{A H}, Asmussen and Hering established a Kesten-Stigum
$L\log L$ type theorem for a class branching diffusion processes
under a condition which is later called a positive regular property
in \cite{A H 1}. In \cite{LRS, LRS2} we established Kesten-Stigum
$L\log L$ type theorems for super-diffusions and branching Hunt
processes respectively.

This paper is a natural continuation of \cite{LRS, LRS2}. The main
purpose of this paper is to establish a strong law of large numbers
for a class of super-diffusions. The main tool of this paper is the
stochastic integral representation of super-diffusions.

Throughout this paper, we will use the following notations. For any
positive integer $k$, $C^k_b(\R^d)$ denotes the family of bounded
functions on $\R^d$ whose partial derivatives of order up to $k$ are
bounded and continuous, $C^k_0(\R^d)$ denotes the family of
functions of compact support on $\R^d$ whose partial derivatives of
order up to $k$ are continuous. For any open set $D\subset \R^d$,
the meanings of $C^k_b(D)$ and $C^k_0(D)$ are similar. We denote by
${\mathcal M}_F(D)$ the space of finite measures on $D$ equipped
with the topology of weak convergence. We will use
$\mathcal{M}_F(D)^0$ to denote the subspace of nontrivial measures
in ${\mathcal M}_F(D)$.  The integral of a function $\varphi$ with
respect to a measure $\mu$ will often be denoted as $\langle
\varphi, \mu\rangle$.

For convenience we use the following convention throughout this
paper: For any probability measure $P$, we also use $P$ to denote
the expectation with respect to $P$.

\subsection{Model}
Suppose that $a_{ij}\in C^1_b(\R^d)$, $i,j=1,\cdots, d$, and that
the matrix $(a_{ij})$ is symmetric and satisfies
$$
0<a|\upsilon|^2\leq\sum_{i,j}a_{ij}\upsilon_i\upsilon_j, \quad
\mbox{for all}\  x\in \R^d  \mbox{ and}\ \upsilon\in\R^d
$$
for some positive constant $a$. We assume that $b_i, i=1, \cdots,
d$, are bounded Borel functions on $\R^d$.
We will use $(\xi, \Pi_x, x\in \R^d)$ to denote a diffusion process
on $\R^d$ corresponding to the operator
$$
L=\frac{1}{2}\nabla\cdot a\nabla+b\cdot \nabla.
$$

In this paper we will always assume  that $D$ is a bounded
domain in $\R^d$. We will use $(\xi^D,\ \Pi_x,\ x\in D)$ to denote
the process obtained by killing $\xi$ upon exiting from $D$, that
is,
$$
\xi^D_t=\left\{\begin{array}{ll} \xi_t & \mbox{ if } t<\tau,\cr
\partial, & \mbox{ if } t\ge \tau,\cr
\end{array}
\right.
$$
where $\tau =\inf\{t>0; \xi_t\notin D\}$ is the first exit time of
$D$ and $\partial$ is a cemetery point. Any function $f$ on $D$ is
automatically extended to $D\cup \{\partial\}$ by setting
$f(\partial)=0$.

We will always assume that $\beta$ is a  bounded Borel function on
$\R^d$. We will use $\{P^D_t\}_{t\ge 0}$ to denote the following
Feynman-Kac semigroup
$$
P^D_tf(x)=\Pi_x\left(\exp\left(\int^t_0\beta(\xi^D_s){\rm
d}s\right)f(\xi^D_t)\right), \quad x\in D.
$$
It is well known that the semigroup $\{P^D_t\}_{t\ge 0}$ is strongly
continuous in $L^2(D)$ and, for any $t>0$, $P^D_t$ has a bounded,
continuous and strictly positive density $p^D(t, x, y)$.

Let $\{\widehat{P}^D_t\}_{t\ge 0}$ be the dual semigroup of
$\{P^D_t\}_{t\ge 0}$ defined by
$$
\widehat{P}^D_tf(x)=\int_Dp^D(t, y, x)f(y){\rm d}y, \quad x\in D.
$$
It is well known that $\{\widehat{P}^D_t\}_{t\ge 0}$ is also
strongly continuous in $L^2(D)$.

Let ${\bf A}$ and $\widehat{{\bf A}}$ be the generators of the
semigroups $\{P^D_t\}_{t\ge 0}$ and $\{\widehat{P}^D_t\}_{t\ge 0}$
in $L^2(D)$ respectively. Let  $\sigma({\bf A})$
($\sigma(\widehat{{\bf A}})$ resp.) denote the spectrum of
 ${\bf A}$ ($\widehat{{\bf A}}$, resp.).  It follows
from Jentzsch's theorem (\cite[Theorem V.6.6, p. 337]{Sc}) and the
strong continuity of $\{P^D_t\}_{t\ge 0}$ and
$\{\widehat{P}^D_t\}_{t\ge 0}$ that the common value $\lambda_1:=
\sup {\rm Re}(\sigma({\bf A}))= \sup {\rm Re}(\sigma(\widehat{{\bf
A}}))$ is an eigenvalue of multiplicity 1 for both ${\bf A}$ and
$\widehat{{\bf A}}$, and that an eigenfunction $\phi$ of ${\bf A}$
associated with $\lambda_1$ can be chosen to be strictly positive
a.\/e. on $D$ and an eigenfunction $\widetilde{\phi}$ of
$\widehat{{\bf A}}$ associated with $\lambda_1$ can be chosen to be
strictly positive a.\/e. on $D$. By \cite[Proposition 2.3]{KS1} we
know that $\phi$ and $\widetilde{\phi}$ are bounded and continuous
on $D$, and they are in fact strictly positive everywhere on $D$. We
choose $\phi$ and $\widetilde{\phi}$ so that
$\int_D\phi(x)\widetilde{\phi}(x){\rm d} x=1$.

Throughout this paper we assume the following

\begin{assumption}\label{assume1}
The semigroups $\{P^D_t\}_{t\ge 0}$ and $\{\widehat{P}^D_t\}_{t\ge
0}$ are intrinsically ultracontractive, that is, for any $t>0$,
there exists a constant $c_t>0$ such that
$$
p^D(t, x, y)\le c_t\phi(x)\widetilde{\phi}(y), \quad \mbox{ for all
} (x, y)\in D\times D.
$$
\end{assumption}

Assumption 1 is a very weak regularity assumption on $D$. It follows
from \cite{KS1, KS2} that Assumption \ref{assume1} is satisfied when
$D$ is a bounded Lipschitz domain.
 For other, more general, examples
of domain $D$ for which Assumption \ref{assume1} is satisfied, we
refer our readers to \cite{KS2} and the references therein.

Define
\begin{equation}\label{def-p-phi}
p^\phi(t, x, y)=\frac{e^{-\lambda_1 t}}{\phi(x)}\ p^D(t, x, y)\
\phi(y).
\end{equation}
Then it follows from \cite[Theorem 2.7]{KS1} that if Assumption
\ref{assume1} holds,  then for any $\sigma>0$ there are positive
constants $C(\sigma)$ and $\nu$ such that
\begin{equation}\label{uniform convergence}
\left|p^\phi(t, x, y)-\phi(y)\widetilde\phi(y)\right|=
\left|\frac{e^{-\lambda_1t}p^D(t, x,y)\phi(y)} {\phi(x)
}-\phi(y)\widetilde\phi(y)\right|\le C(\sigma) e^{-\nu
t}\phi(y)\widetilde\phi(y),\quad x,y\in D, t>\sigma.
\end{equation}
By the definition of $\phi$ and $\widetilde{\phi}$, it is easy to
check that, for any $t>0$, $p^\phi(t, \cdot, \cdot)$ is a
probability density and that $\phi\widetilde\phi$ is its unique
invariant probability density. \eqref{uniform convergence} shows
that $p^\phi(t,\cdot,x)$ converges to $\phi(x)\widetilde\phi(x)$
uniformly with exponential rate. Denote by $P_t^\phi$ the semigroup
with density  $ p^\phi(t, \cdot, \cdot)$ and $\Pi_h^\phi$ the
 probability generated by $(P^\phi_t)_{t\geq 0}$ with initial
distribution $h(x){\rm d}x$ on $D$. Then $(\xi^D,
\Pi^\phi_{\phi\widetilde\phi})$ is a diffusion with initial
distribution $\phi(x)\widetilde\phi(x){\rm d}x$.

The super-diffusion $(X, \mathbb P_\mu), \mu\in \mathcal M_F(D)^0$,
we are going to study is a $(\xi^D,\ \psi(\lambda) -
\beta\lambda)$-super-process, which is a measure-valued Markov
process with underlying spatial motion $\xi^D$, branching rate ${\rm
d}t$ and branching mechanism $\psi(\lambda)-\beta\lambda$, where
$$
\psi(x, \lambda)= \frac{1}{2}\alpha(x)\lambda^2+
\int_0^\infty\big(e^{-r\lambda}-1+\lambda r\big)n(x, {\rm
d}r),\qquad \lambda>0,
$$
for some nonnegative bounded measurable function $\alpha$ on $D$ and
for some $\sigma$-finite kernel $n$ from $(D,\ {\cal B}(D))$ to $(
\R_+, {\cal B}(\R_+))$, that is, $n(x, {\rm d}r)$ is a
$\sigma$-finite measure on $\R_+$ for each fixed $x$, and $n(\cdot,
B)$ is a measurable function for each Borel set $B\subset \R_{+}$.
The measure $\mu$ here is the initial value of $X$.  In this paper
we will always assume that
\begin{equation}\label{assumption on bm}
\sup_{x\in D}\int_0^\infty (r\wedge r^2)n(x,{\rm d}r)<\infty.
\end{equation}
Note that this assumption implies, for any fixed $\lambda>0$,
$\psi(\cdot,\lambda)$ is bounded on $D$. Define a new kernel
$n^\phi(x, {\rm d}r)$ from $(D,\ {\cal B}(D))$ to $(\R_+, {\cal
B}(\R_+))$ such that for any nonnegative measurable function $f$ on
$\R_+$,
\begin{equation}\label{phi-change}
\int_0^\infty f(r)n^\phi(x,{\rm d}r)=\int_0^\infty f(r\phi(x))n(x,
{\rm d}r), \quad x\in D.
\end{equation}
Then, by \eqref{assumption on bm} and the boundedness of
$\phi$, $n^\phi$ satisfies
\begin{equation}\label{assumption on transformed bm}
\sup_{x\in D}\int_0^\infty (r\wedge r^2)n^\phi(x,{\rm d}r)<\infty.
\end{equation}

\subsection{Stochastic Integral Representation and Main Result}

Let $(\Omega,\ \mathcal{F}, \mathbb P_\mu, \mu\in
\mathcal{M}_{F}(D)^0)$ be the underlying probability space equipped
with the filtration $(\mathcal{F}_t)$, which is generated by $X$ and
is completed as usual with the $\mathcal{F}_\infty-$measurable and
$\mathbb P_\mu-$negligible sets for every $\mu\in
\mathcal{M}_{F}(D)^0$.  It is known (cf. \cite[Section 6.1]{D}) that
the super-diffusion $X$ is a
solution to the following martingale problem: for any $\varphi\in
C_0^2(D)$ and $h\in C^2_b(\mathbb{R}),$
\begin{eqnarray}\label{martingale problem}
h(\langle \varphi, X_t\rangle)-h(\langle\varphi,\mu\rangle)
-\int_0^th'(\langle\varphi, X_s \rangle)\langle
{\bf A}\varphi, X_s\rangle {\rm d}s-\frac{1}{2}\int_0^th''(
\langle\varphi, X_s \rangle)\langle\alpha\varphi^2, X_s\rangle{\rm d}s\\
-\int_0^t\int_D\int_{(0, \infty)}\big(h(\langle \varphi,
X_s\rangle+r\varphi(x))-h(\langle \varphi, X_s\rangle)-h'(\langle
\varphi, X_s\rangle)r\varphi(x)\big)n(x, {\rm d}r) X_s({\rm d}x){\rm
d}s\nonumber
\end{eqnarray}
is a martingale. Let $J$ denote the set of all jump times of $X$ and
$\delta$ denote the Dirac measure.  It is easy to see from the
martingale problem \eqref{martingale problem} that the only possible
jumps of $X$ are point measures $r\delta_{x}(\cdot)$ with  $r$ being
a positive real number and $x\in D$. \eqref{martingale problem}
implies that the compensator of the random measure
$$
   N:=\sum_{s\in  J}\delta_{(s, \triangle X_s)}
$$
is a random measure $\widehat{N}$ on $\R_+\times\mathcal{M}_F(D)^0$
such that for any nonnegative predictable function $F$ on
$\R_+\times \Omega\times \mathcal{M}_F(D)^0,$
\begin{eqnarray}\label{poisson point definition}
\int\int F(s,\omega, \nu)\widehat{N}({\rm d}s,{\rm d}\nu)
=\int_0^\infty {\rm d}s\int_D X_s({\rm d}x) \int_0^\infty
F(s,\omega, r\delta_x)n(x,{\rm d}r),
\end{eqnarray}
where $n(x,{\rm d}r)$ is the kernel of the branching mechanism
$\psi$.  Therefore we have
\begin{eqnarray}\label{expectation identity}
\mathbb P_\mu\left[\sum_{s\in J}F(s,\omega,\triangle X_s)\right]
=\mathbb P_\mu\int_0^\infty {\rm d}s\int_D X_s({\rm d}x)
\int_0^\infty  F(s,\omega, r\delta_x)n(x,{\rm d}r).
\end{eqnarray}
See \cite[p. 111]{D}.   Let $F$ be a predictable function on
$\mathbb{R}_+\times\Omega\times\mathcal{M}_F(D)^0$ satisfying
$$
\mathbb P_\mu\bigg[\bigg(\sum_{s\in[0,t],s\in J}F(s, \Delta
X_s)^2\bigg)^{1/2}\bigg]<\infty, \qquad \mbox{ for all } \mu\in
{\mathcal M}_F(D)^0.
$$
Then the stochastic integral of $F$ with respect to the compensated
random measure $N-\widehat{N}$
$$
\int_0^t\int_{\mathcal{M}_F(D)^0}F(s, \nu)(N-\widehat{N})({\rm d}s,
{\rm d}\nu),
$$
can be defined (cf. \cite{LM} and the reference therein) as the
unique purely discontinuous martingale (vanishing at time $0$) whose
jumps are indistinguishable from $1_J(s)F(s,\Delta X_s).$

Suppose that $\varphi$ is a measurable function on $\R_+\times D$.
Define
\begin{eqnarray}F_\varphi(s,\nu):=\int_D\varphi(s,x)\nu({\rm
d}x),\quad \nu\in \mathcal{M}_F(D)
\end{eqnarray}
whenever the integral above makes sense.  We write
\begin{eqnarray}\label{martingale definition}
S^J_t(\varphi)=\int_0^t\int_D\varphi(s,x)S^J({\rm d}s,{\rm d}x)
:=\int_0^t\int_{\mathcal{M}_F(D)^0}F_\varphi(s,\nu)(N-\widehat{N})({\rm
d}s,{\rm d}\nu),
\end{eqnarray}
whenever the right hand of $(\ref{martingale definition})$ makes
sense. If $\varphi$ is bounded on $\R_+\times D$, then
$S^J_t(\varphi)$ is well defined.  Indeed, we only need to check that
\begin{eqnarray}\label{int-finit}
\mathbb P_\mu\left[\left(\sum_{s\in[0,t],s\in J}F_\varphi(s, \Delta
X_s)^2\right)^{1/2}\right]<\infty, \qquad \mbox{ for all } \mu\in
{\mathcal M}_F(D)^0.
\end{eqnarray}
Note that, for any $\mu\in {\mathcal M}_F(D)^0$,
\begin{equation*}
\begin{array}{rl}
&\mathbb P_\mu\left[\left(\sum_{s\in[0,t],s\in J}
F_\varphi(s, \Delta X_s)^2\right)^{1/2}\right]\\
=&\mathbb P_\mu\left[\left(\sum_{s\in[0,t],s\in J}
\left(\int\varphi(s,x)(\Delta X_s)({\rm d}x)\right)^2\right)^{1/2}\right]\\
\le&\displaystyle\|\varphi\|_\infty \mathbb P_\mu\left[\left(
\sum_{s\in[0,t],s\in J}\langle 1, \Delta X_s\rangle^2 I_{\{\langle
1, \Delta X_s\rangle\le 1\}}+\sum_{s\in[0,t],s\in J} \langle 1,
\Delta X_s\rangle^2 I_{\{\langle 1, \Delta X_s\rangle> 1\}}
\right)^{1/2}\right]\\
\le&\displaystyle\|\varphi\|_\infty \mathbb P_\mu\left[\left(
\sum_{s\in[0,t],s\in J}\langle 1, \Delta X_s\rangle^2 I_{\{
\langle 1, \Delta X_s\rangle\le 1\}}\right)^{1/2}\right]\\
&+\displaystyle\|\varphi\|_\infty \mathbb P_\mu\left[\left(
\sum_{s\in[0,t],s\in J}\langle 1,\ \Delta X_s\rangle^2 I_{\{ \langle
1, \Delta X_s\rangle> 1\}}\right)^{1/2}\right].
\end{array}
\end{equation*}
Here and throughout this paper, for any set $A$, $I_A$ stands for the
indicator function of $A$.
Using the first two displays on \cite[p. 203]{LM}, we get
\eqref{int-finit}. Thus for any bounded function $\varphi$ on
$\mathbb R_+\times D$. $(S^J_t(\varphi))_{t\geq 0}$ is a martingale.

For any $\varphi\in
 C^2_0(D)$ and $\mu\in {\mathcal M}_F(D)$,
\begin{equation}\label{Martingale-expres1}
\langle \varphi, X_t\rangle = \langle\varphi, \mu\rangle
+S^J_t(\varphi)+ S^C_t(\varphi) + \int_0^t\langle{\bf
A}\varphi, X_s\rangle {\rm d}s,
\end{equation}
where $S^C_t(\varphi)$ is a continuous local martingale with
quadratic variation
\begin{equation}\label{quadratic variation}
\langle S^C(\varphi)\rangle_t=\int_0^t\langle \alpha\varphi^2, X_s\rangle {\rm d}s.
\end{equation}
In fact, according to \cite{F1, F2}, the above is still valid when
${\bf A}$ is replaced by ${\bf L}+\beta$, where ${\bf L}$ is the
weak generator of $\xi^D$ in the sense of \cite[Section 4]{F1}.
Using this, \cite[Corollary 2.18]{F2} and applying a limit argument,
one can show that for any bounded function $g$ on $D$,
\begin{equation}\label{dual representation}
\langle g, X_t\rangle=\langle P^D_tg,\mu\rangle+\int_0^t\int_D
P^D_{t-s}g(x)S^J({\rm d}s, {\rm d}x)+\int_0^t\int_D
P^D_{t-s}g(x)S^C({\rm d}s, {\rm d}x).
\end{equation}
In particular, taking $g=\phi$ in \eqref{dual representation}, where
$\phi$ is the positive eigenfunction of $\bf A$ defined in Section
1.1, we get that
\begin{eqnarray}\label{stochatic r for mart}
e^{-\lambda_1 t}\langle \phi, X_t\rangle =\langle \phi,
\mu\rangle+\int_0^t e^{-\lambda_1 s}\int_D \phi(x)S^J({\rm d} s,
{\rm d}x)+\int_0^t e^{-\lambda_1 s}\int_D \phi(x)S^C({\rm d} s, {\rm
d}x).
\end{eqnarray}
Set $M_t(\phi):=e^{-\lambda_1 t}\langle\phi, X_t\rangle$.  Then
$M_t(\phi), t\ge 0,$ is a nonnegative martingale.  Denote by
$M_\infty(\phi)$ the almost sure limit of $M_t(\phi)$ as
 $t\to\infty$. In \cite{LRS}, we studied the relationship
between the degeneracy property of $M_\infty(\phi)$ and the function
$l$:
\begin{equation}\label{m}
l(y):=\int_1^\infty r\ln r n^\phi(y, {\rm d}r).
\end{equation}

\begin{thm}\label{degeneracy theorem}\cite[Theorem 1.1]{LRS}
Suppose that Assumption \ref{assume1} holds, $\lambda_1>0$ and that
$X$ is a $(\xi^D, \psi(\lambda)-\beta\lambda)-$super-diffusion. Then
the following assertions hold:
\begin{description}
\item{(1)} If $\int_Dl(y)\widetilde{\phi}(y){\rm d}y<\infty$, then
$M_\infty(\phi)$ is non-degenerate under $\mathbb P_\mu$ for any
$\mu\in \mathcal M_F(D)^0$, and $M_\infty(\phi)$ is also the
$L^1(\mathbb P_\mu)$ limit of $M_t(\phi)$.
\item{(2)} If
$\int_Dl(y)\widetilde{\phi}(y){\rm d}y=\infty$, then
$M_\infty(\phi)=0,\, \mathbb P_\mu{\rm -a.s.}$ for any
$\mu\in\mathcal M_F(D)^0$.
\end{description}
\end{thm}

\begin{remark} In \cite[Theorem 1.1]{LRS}, we only stated that
in case (1) under the extra assumption $\alpha\equiv 0$, $M_\infty(\phi)$ is
non-degenerate under $\mathbb P_\mu$ for any $\mu\in \mathcal
M_F(D)^0$. But actually in this case we have $\mathbb P_\mu
M_\infty(\phi)=\mathbb P_\mu M_0(\phi)$ (see \cite[Lemma 3.4]{LRS}),
and therefore $M_t(\phi)$ converges to $M_\infty(\phi)$ in
$L^1(\mathbb P_\mu)$.

For general $\alpha\ge 0$, by the $L^2$ maximum inequality, and
using the fact that $\alpha$ and $\phi$ are bounded in $D$, we have
\begin{equation}
\begin{split}
&\mathbb P_\mu\left[\sup_{t\geq 0}\left(\int_0^t e^{-\lambda_1 s}\int_D
\phi(x)S^C({\rm d} s, {\rm d}x)\right)^2\right]\\
\leq& 4\sup_{t\geq 0}\mathbb P_\mu \left(\int_0^t e^{-\lambda_1 s}\int_D
\phi(x)S^C({\rm d} s, {\rm d}x)\right)^2\\
=& 4 \mathbb P_\mu \int_0^\infty e^{-2\lambda_1 s}{\rm d}s\int_D \alpha(x)
\phi^2(x)X_s({\rm d}x)\\
=& 4 \int_0^\infty e^{-\lambda_1 s}{\rm d}s\int_D\phi(y)\mu({\rm d}y)
\int_Dp^\phi(s,y,x)\alpha(x)\phi(x){\rm d}x\\
<&\infty.
\end{split}
\end{equation}
Thus the martingale $\left(\int_0^t e^{-\lambda_1 s}\int_D
\phi(x)S^C({\rm d} s, {\rm d}x)\right)_{t\geq 0}$ converges almost
surely and in $L^1(\mathbb P_\mu)$.  Denote the limit by
$\int_0^\infty e^{-\lambda_1 s}\int_D \phi(x)S^C({\rm d} s, {\rm
d}x)$.  Furthermore, we obtain that when $\lambda_1>0$ and
$\int_Dl(x)\widetilde\phi(x){\rm d}x<\infty$, the martingale
$\int_0^t e^{-\lambda_1 s}\int_D \phi(x)S^J({\rm d}s, {\rm d}x)$
converges almost surely and in $L^1(\mathbb P_\mu)$ as well.  Denote
the limit by $\int_0^\infty e^{-\lambda_1 s}\int_D \phi(x)S^J({\rm
d}s, {\rm d}x)$.
Thus it follows from \eqref{stochatic r for mart} that $M_t(\phi)$ converges to
a non-degenerate $M_\infty(\phi)$ $\mathbb P_\mu$-almost surely and
in $L^1(\mathbb P_\mu)$ for every $\mu \in\mathcal M_F(D)^0$.
\end{remark}

The main goal of this paper is to establish the following almost
sure convergence result.

\begin{thm}\label{new main theorem}
Suppose that Assumption \ref{assume1} holds, $\lambda_1>0$ and that
$X$ is a $(\xi^D, \psi(\lambda)-\beta\lambda)-$super-diffusion. Then
there exists $\Omega_0\subset\Omega$ of probability one (that is,
$\mathbb P_\mu (\Omega_0)=1$ for every $\mu\in \mathcal M_F(D)^0$)
such that, for every $\omega\in\Omega_0$ and for every nontrivial
nonnegative bounded Borel function $f$ on $D$ with compact support
whose set of discontinuous points has zero Lebesgue measure, we have
\begin{eqnarray}\label{new main result}
\lim_{t\rightarrow\infty}\frac{\langle f,
X_t\rangle(\omega)}{\mathbb P_\mu\langle f, X_t\rangle} =
\frac{M_\infty(\phi)(\omega)}{\langle\phi,\mu\rangle}.
\end{eqnarray}
\end{thm}

As a consequence of this theorem we immediately get the following

\begin{corollary}
Suppose that Assumption \ref{assume1} holds, $\lambda_1>0$ and that
$X$ is a $(\xi^D, \psi(\lambda)-\beta\lambda)-$super-diffusion. Then
there exists $\Omega_0\subset\Omega$ of probability one (that is,
$\mathbb P_\mu (\Omega_0)=1$ for every $\mu\in \mathcal M_F(D)^0$)
such that, for every $\omega\in \Omega_0$ and every relatively
compact Borel subset $B$ in $D$ of positive Lebesgue measure whose
boundary is of Lebesgue measure zero, we  have
$$
\lim_{t\to\infty}\frac{X_t(B)(\omega)}{\mathbb P_\mu[X_t(B)]}=
\frac{M_\infty(\phi)(\omega)}{\langle\phi,\mu\rangle}.
$$
\end{corollary}

\begin{remark}
(i) Although we assumed in this paper that the underlying motion
$\xi^D$ is a diffusion process in a bounded domain $D$, the
arguments of this paper can be easily extended to the case when the
underlying
motion is a Hunt process on a locally compact
separable metric space $E$ satisfying \cite[Assumption 1.1]{LRS2}
for some measure $m$ with full support and with $m(E)<\infty$, and
the analogue of Assumption 1 above.

(ii) In \cite{CRW, E, EHK, EW, W}, the branching mechanism is
assumed to be binary, while in the present paper we deal with
general branching mechanism. \cite{ET} considers a general branching
mechanism under a $(1+\theta)$-moment condition, $\theta>0$,
while in the present paper, we only assume a
$L\log L$
condition. In
\cite{CRW} the underlying motion is assumed to be a symmetric Hunt
process, while in the present paper, our underlying process needs
not be symmetric.
\end{remark}

\section{Proof of Theorem \ref{new main theorem}}

A main step in  proving  Theorem \ref{new main theorem} is the
following result.

\begin{thm}\label{main theorem}
Suppose that Assumption \ref{assume1} holds, $\lambda_1>0$ and that
$X$ is a $(\xi^D, \psi(\lambda)-\beta\lambda)-$super-diffusion. Then
for any $f\in\mathcal B_b(D)$ and $\mu\in \mathcal M_F(D)^0$,
\begin{eqnarray}\label{main result}
\lim_{t\rightarrow\infty}e^{-\lambda_1 t}\langle\phi f, X_t\rangle =
M_\infty(\phi)\int_D\phi(y)\widetilde\phi(y)f(y){\rm d}y,\qquad
\mathbb P_\mu{\rm -a.s.}
\end{eqnarray}
\end{thm}

We will prove this result first. According to  Theorem
\ref{degeneracy theorem} (2), when $\int_D\widetilde\phi(x)l(x){\rm
d}x=\infty$, we have
$$
M_\infty(\phi)=\lim_{t\rightarrow\infty}e^{-\lambda_1 t}\langle\phi,
X_t\rangle=0,\qquad \mathbb P_\mu{\rm -a.s.}
$$
For any $f\in \mathcal B^+_b(D)$,
$$
\limsup_{t\rightarrow\infty}e^{-\lambda_1 t}\langle\phi f,
X_t\rangle\leq \|f\|_{\infty} \limsup_{t\rightarrow\infty}
e^{-\lambda_1 t}\langle\phi, X_t\rangle=0
$$
and \eqref{main result} follows immediately from the nonnegativity
of $f$.

It remains to prove the case when $\int_D\widetilde\phi(x)l(x){\rm
d}x<\infty$.  In the remainder of this section, we assume that the
assumptions of Theorem \ref{main theorem} hold and that
$f\in\mathcal B_b^+(D)$ is fixed. Define
$$
S({\rm d}s,{\rm d}x)=S^J({\rm d}s,{\rm d}x)+S^C({\rm d}s,{\rm d}x).
$$
Using \eqref{dual representation}, we get
\begin{equation}\label{equation 1}
e^{-\lambda_1 t}\langle\phi f, X_t\rangle=e^{-\lambda_1 t}\langle
P^D_t(\phi f),\mu\rangle+e^{-\lambda_1 t}
\int_0^t\int_D(P^D_{t-s}(\phi f))(x)S({\rm d}s,{\rm d}x).
\end{equation}
It follows from \eqref{uniform convergence} that
$$
\lim_{t\rightarrow\infty}e^{-\lambda_1 t}\langle P^D_t(\phi
f),\mu\rangle=\langle\phi,\mu\rangle\int_D\phi(x)\widetilde\phi(x)f(x){\rm
d}x.
$$
Therefore, to prove Theorem \ref{main theorem}, it suffices to show
that
\begin{equation}\label{leftover}
\lim_{t\rightarrow\infty}e^{-\lambda_1 t}
\int_0^t\int_D(P^D_{t-s}(\phi f))(x)S({\rm d}s,{\rm
d}x)=\int_D\phi\widetilde\phi(x)f(x){\rm d}x\int_0^\infty
e^{-\lambda_1 s}\int_D\phi(x) S({\rm d}s, {\rm d}x),\qquad \mathbb
P_\mu{\rm -a.s.}
\end{equation}

To prove the above result, we need some lemmas first.

\begin{lemma}\label{finite lemma}
For any $\mu\in \mathcal M_F(D)^0$, we have
\begin{equation}\label{finite 1}
\sum_{s>0}I_{\{\Delta X_s(\phi)>e^{\lambda_1 s}\}} < \infty ,\qquad
\mathbb P_\mu{\rm -a.s.}
\end{equation}
\end{lemma}

{\bf Proof:}  First note that \eqref{assumption on transformed bm}
implies  that $\sup_{x\in D}\int^\infty_1 rn^\phi(x, {\rm
d}r)<\infty$. It is well known that for any $g\in\mathcal B^+(D)$,
\begin{equation}\label{mean-formu}
\mathbb P_\mu\langle g, X_t\rangle=\langle P^D_tg,\mu\rangle.
\end{equation}
By the definition of $n^\phi$ given by
\eqref{phi-change}, we have
\begin{eqnarray*}
&&\mathbb P_\mu \left[\sum_{s>0}I_{\{\Delta X_s(\phi)>e^{\lambda_1 s}\}}\right]\\
&=&\mathbb P_\mu\int_0^\infty {\rm d}s\int_D\int_0^\infty I_{\{r\phi(x) >
e^{\lambda_1 s}\}}X_s({\rm d}x)n(x,{\rm d}r)\\
&=& \int_0^\infty ds\int_D \mu({\rm d}y)\int_Dp^D(s, y, x){\rm
d}x\int_{e^{\lambda_1 s}}^\infty n^\phi(x, {\rm d}r).
\end{eqnarray*}
It follows from \eqref{uniform convergence} that there is a
constant $C>0$ such that
\begin{equation}\label{uniform boundary}
p^D(t,y,x)\leq Ce^{\lambda_1 t}\phi(y)\widetilde\phi(x),\qquad
\forall\ t>1,\ x,y\in D.
\end{equation}
Therefore we have,
\begin{eqnarray*}
&&\mathbb P_\mu \left[\sum_{s>0}I_{\{\Delta X_s(\phi) >
e^{\lambda_1 s}\}}\right]\\
&=&\int_0^1 {\rm d}s\int_D \mu({\rm d}y)\int_Dp^D(s, y, x){\rm d}x
\int_{e^{\lambda_1 s}}^\infty n^\phi(x, {\rm d}r)+ \int_1^\infty
{\rm d}s\int_D \mu({\rm d}y)\int_Dp^D(s, y, x)
{\rm d}x\int_{e^{\lambda_1 s}}^\infty n^\phi(x, {\rm d}r)\\
&\leq &e^{\|\beta\|_\infty}\mu(D)\left\|\int_1^\infty n^\phi
(x, {\rm d}r)\right\|_\infty+C\int_0^\infty e^{\lambda_1 s}
{\rm d}s\int_D\phi(y)\mu({\rm d}y)\int_D\widetilde\phi(x)
{\rm d}x\int_{e^{\lambda_1 s}}^\infty n^\phi(x, {\rm d}r)\\
&=:&C_1+C\langle\phi,\mu\rangle\int_D\widetilde\phi(x){\rm d}x
\int_0^\infty e^{\lambda_1 s}{\rm d}s\int_{e^{\lambda_1 s}}^\infty
n^\phi(x, {\rm d}r).
\end{eqnarray*}
Using Fubini's theorem we get,
$$
\int_0^\infty e^{\lambda_1 s}{\rm d}s \int_{e^{\lambda_1 s}}^\infty
n^\phi(x, {\rm d}r)\le\frac{1}{\lambda_1}\int^\infty_1 rn^\phi(x,
{\rm d}r).
$$
Hence we have
\begin{eqnarray*}
\mathbb P_\mu \Big[\sum_{s>0}I_{\{\Delta X_s(\phi)> e^{\lambda_1
s}\}}\Big]\leq \frac{C\langle\phi,\mu\rangle}{\lambda_1}
\int_D\widetilde\phi(x){\rm d}x\int_1^\infty rn^{\phi}(x, {\rm d}r)
<\infty.
\end{eqnarray*}
Consequently, \eqref{finite 1} holds.\qed

Define
$$
N_\phi^{(1)}:=\sum_{0<\triangle X_s(\phi)<e^{\lambda_1s}}\delta_{(s,
\,\, \triangle X_s)}\quad \mbox{and}\quad
N_\phi^{(2)}:=\sum_{\triangle X_s(\phi)\ge
e^{\lambda_1s}}\delta_{(s, \,\, \triangle X_s)},
$$
and denote the compensators of $N_\phi^{(1)}$  and $N_\phi^{(2)}$ by
$\widehat N_\phi^{(1)}$ and $\widehat N_\phi^{(1)}$ respectively.
Then for any nonnegative predictable function $F$ on
$\mathbb{R}_+\times \Omega\times \mathcal{M}_F(D)$,
\begin{equation}\label{identity1}
\int_0^\infty \int F(s, \nu)\widehat N_\phi^{(1)}({\rm d}s, {\rm
d}\nu) =\int_0^\infty {\rm d}s\int_DX_s({\rm
d}x)\int_0^{e^{\lambda_1s}} F(s, r\phi(x)^{-1}\delta_x)n^{\phi}(x,
{\rm d}r),
\end{equation}
and
\begin{equation}\label{identity1(2)}
\int_0^\infty \int F(s, \nu)\widehat N_\phi^{(2)}({\rm d}s, {\rm
d}\nu) =\int_0^\infty {\rm d}s\int_DX_s({\rm
d}x)\int_{e^{\lambda_1s}}^\infty F(s,
r\phi(x)^{-1}\delta_x)n^{\phi}(x, {\rm d}r).
\end{equation}
Let $J^{(1)}_\phi$ denote the set of jump times of $ N_\phi^{(1)}$,
and $J^{(2)}_\phi$ the set of jump times of $ N_\phi^{(2)}$. Then
\begin{equation}\label{identity1a}
\int_0^\infty \int F(s, \nu) N_\phi^{(1)}({\rm d}s, {\rm d}\nu)
=\sum_{s\in J^{(1)}_\phi}F(s,\omega,\triangle
 X_s),
\end{equation}
\begin{equation}\label{identity1(2)a}
\int_0^\infty \int F(s, \nu) N_\phi^{(2)}({\rm d}s, {\rm d}\nu)
=\sum_{s\in J^{(2)}_\phi}F(s,\omega,\triangle
 X_s),
\end{equation}
\begin{eqnarray}\label{expectation identity2}
\mathbb P_\mu\left[\sum_{s\in J^{(1)}_\phi}F(s,\omega,\triangle
X_s)\right] =\mathbb P_\mu\int_0^\infty {\rm d}s\int_D X_s({\rm d}x)
\int_0^{e^{\lambda_1s}}   F(s,\omega,
r\phi(x)^{-1}\delta_x)n^\phi(x,{\rm d}r),
\end{eqnarray}
and
\begin{eqnarray}\label{expectation identity2(2)}
\mathbb P_\mu\left[\sum_{s\in J^{(2)}_\phi}F(s,\omega,\triangle
X_s)\right] =\mathbb P_\mu\int_0^\infty {\rm d}s\int_D X_s({\rm d}x)
\int_{e^{\lambda_1s}}^\infty   F(s,\omega,
r\phi(x)^{-1}\delta_x)n^\phi(x,{\rm d}r).
\end{eqnarray}
We can construct two martingale measures
$S^{J,(1)}({\rm d}s, {\rm d}x)$ and $S^{J,(2)}({\rm d}s, {\rm d}x)$
respectively from $N^{(1)}_\phi({\rm d}s, {\rm d}\nu)$ and
$N^{(2)}_\phi({\rm d}s, {\rm d}\nu)$, similar to the way we
 constructed $S^J({\rm d}s, {\rm d}x)$ from $N({\rm d}s, {\rm d}\nu)$.
Then for any bounded measurable
function $g$ on $\R_+\times D$,
\begin{equation}\label{martingale-int2}
S^{J, (1)}_t(g)=\int^t_0\int_Dg(s,x) S^{(1)}({\rm d}s, {\rm d}x)=
\int^t_0\int_{{\cal M}_F(D)}F_g(s,\nu)(N^{(1)}_\phi-\widehat
N^{(1)}_\phi)({\rm d}s, {\rm d}\nu),
\end{equation}
and
\begin{equation}\label{martingale-int2(2)}
S^{J, (2)}_t(g)=\int^t_0\int_Dg(s,x) S^{(2)}({\rm d}s, {\rm d}x)
=\int^t_0\int_{{\cal M}_F(D)}F_g(s,\nu)(N^{(2)}_\phi-\widehat
N^{(2)}_\phi)({\rm d}s, {\rm d}\nu),
\end{equation}
where $F_g(s,\nu)=\int g(s,x)\nu({\rm d}x).$

For any $m,n\in\mathbb N$,  $\sigma>0$  and $f\in\mathcal B_b^+(D)$,
define
$$
H_{(n+m)\sigma}(f):=e^{-\lambda_1 (n+m)\sigma}
\int_0^{(n+m)\sigma}\int_DP^D_{(n+m)\sigma-s}(\phi f)(x)
S^{J, (1)}({\rm d}s,{\rm d}x)
$$
and
$$
L_{(n+m)\sigma}(f):= e^{-\lambda_1 (n+m)\sigma}
\int_0^{(n+m)\sigma}\int_D P^D_{(n+m)\sigma-s}(\phi f)(x)
S^{J, (2)}({\rm d}s, {\rm d}x).
$$

\begin{lemma}\label{limit lemma 1}
If $\int_Dl(x)\widetilde\phi(x){\rm d}x<\infty$, then for any $m
\in\mathbb N$, $ \sigma>0$, $\mu\in {\mathcal M}_F(D)^0$ and
$f\in\mathcal B_b^+(D)$,
\begin{equation}\label{finite 2}
\sum_{n=1}^\infty \mathbb P_\mu\left[H_{(n+m)\sigma}(f) -\mathbb
P_\mu(H_{(n+m)\sigma}(f)\big|\mathcal F_{n\sigma})\right]^2<\infty
\end{equation}
and
\begin{equation}\label{limit1}
\lim_{n\rightarrow\infty}H_{(n+m)\sigma}(f)- \mathbb
P_\mu[H_{(n+m)\sigma}(f)\big|\mathcal F_{n\sigma}]=0,\qquad \mathbb
P_\mu {\rm-a.s.}
\end{equation}
\end{lemma}

{\bf Proof:} Since $P^D_t(\phi f)$ is  bounded in $[0,T]\times D$
for any $T>0$, the process
$$
H_{t}(f):=e^{-\lambda_1 (n+m)\sigma} \int_0^{t} \int_D
P^D_{(n+m)\sigma-s}(\phi f)(x) S^{J, (1)}({\rm d}s, {\rm d}x),\qquad
t\in [0, (n+m)\sigma]
$$
is a martingale with respect to $(\mathcal F_{t})_{t\leq
(n+m)\sigma}$. Thus
$$
\mathbb P_\mu(H_{(n+m)\sigma}(f)\big|\mathcal
F_{n\sigma})=e^{-\lambda_1 (n+m)\sigma} \int_0^{n\sigma}\int_D
P^D_{(n+m)\sigma-s}(\phi f)(x) S^{J, (1)}({\rm d}s, {\rm d}x),
$$
and hence
$$
H_{(n+m)\sigma}(f)-\mathbb P_\mu(H_{(n+m)\sigma}(f) \big|\mathcal
F_{n\sigma}) =e^{-\lambda_1 (n+m)\sigma}\int_{n\sigma}^{(n+m)
\sigma}\int_D P^D_{(n+m)\sigma-s}(\phi f)(x) S^{J, (1)}({\rm
d}s,{\rm d}x).
$$
Since
\begin{eqnarray*}
M_t&:=&e^{-\lambda_1 (n+m)\sigma}\int_{n\sigma}^t\int_DP^D_{(n+m)
\sigma-s}(\phi f)(x) S^{J, (1)}({\rm d}s,{\rm d}x)\\
&=&\int_{n\sigma}^t F_{e^{-\lambda_1
(n+m)\sigma}P^D_{(n+m)\sigma-\cdot}(\phi f)}(s, \nu) (N^{(1)}_\phi-
\widehat{N}^{(1)}_\phi)({\rm d}s, {\rm d}\nu), \quad t\in[n\sigma, (n+m)\sigma]
\end{eqnarray*}
is a martingale with quadratic variation
\begin{eqnarray*}
\int_{n\sigma}^t F_{e^{-\lambda_1
(n+m)\sigma}P^D_{(n+m)\sigma-\cdot}(\phi f)}(s, \nu)^2
\widehat{N}^{(1)}_\phi({\rm d}s, {\rm d}\nu),
\end{eqnarray*}
we have
\begin{eqnarray}\label{inequality 1}
&&\mathbb P_\mu\left[H_{(n+m)\sigma}(f)-\mathbb
P_\mu(H_{(n+m)\sigma}(f)\big|\mathcal
F_{n\sigma})\right]^2\nonumber\\
&=& \mathbb P_\mu\int_{n\sigma}^{(n+m)\sigma}F_{e^{-\lambda_1
(n+m)\sigma}P^D_{(n+m)\sigma-\cdot}(\phi f)}(s, \nu)^2
\widehat{N}^{(1)}_\phi({\rm d}s,\ {\rm d}\nu)\nonumber\\
&=&\mathbb P_\mu\bigg[\sum_{s\in \tilde J^{(1)}_{n,m}}
F_{e^{-\lambda_1 (n+m)\sigma}P^D_{(n+m)\sigma-\cdot}(\phi f)}(s,
\Delta X_s)^2\bigg],
\end{eqnarray}
where $\tilde J^{(1)}_{n,m}= J^{(1)}_\phi\bigcap\ [n\sigma,
(n+m)\sigma]$.  Note that for any $f\in\mathcal B_b(D)$,
$\|P^\phi_tf\|_{\infty}\le \|f\|_{\infty}$ for all $t\ge 0$, which
is equivalent to
\begin{equation}\label{uniform boundary2}
P^D_t(\phi f)(y)\leq \|f\|_{\infty}e^{\lambda_1 t}\phi(y), \quad
\forall\ t\ge 0,\  y\in D.
\end{equation}
Using \eqref{identity1} and \eqref{expectation identity2}, we obtain
\begin{eqnarray*}
&&\mathbb P_\mu\bigg[\sum_{s\in \widetilde J_{n,m}}F_{e^{-\lambda_1
(n+m)\sigma} P^D_{(n+m)\sigma-\cdot}(\phi f)}(s, \Delta X_s)^2\bigg]\\
&=&\mathbb P_{\mu}\int_{n\sigma}^{(n+m)\sigma} {\rm d}s\int_D
X_s({\rm d}x)\int_0^{e^{\lambda_1 s}}F^2_{e^{-\lambda_1 (n+m)
\sigma}P^D_{(n+m)\sigma-\cdot}(\phi f)}(s, r
\phi(x)^{-1}\delta_x)n^\phi(x,{\rm d}r)\\
&=&e^{-2\lambda_1 (n+m)\sigma}\int_{n\sigma}^{(n+m)\sigma} {\rm d}s
\int_D \mu({\rm d}y)\int_D p^D(s, y, x){\rm d}x\int_0^{e^{\lambda_1
s}} [P^D_{(n+m)\sigma-s} (\phi f)(x)\phi(x)^{-1}]^2r^2n^\phi(x,{\rm d}r)\\
&\le& \|f\|_\infty^2 \int_{n\sigma}^{(n+m)\sigma}
e^{-2\lambda_1s}{\rm d}s \int_D \mu({\rm d}y)\int_D p^D(s, y, x){\rm
d}x\int_0^{e^{\lambda_1 s}} r^2n^\phi(x,{\rm d}r),
\end{eqnarray*}
where in the second equality we used the fact that
\begin{eqnarray}\label{formu-F}
F_{e^{-\lambda_1 (n+m)\sigma}P^D_{(n+m)\sigma-\cdot}(\phi f)}(s,
r\phi(x)^{-1}\delta_x)= re^{-\lambda_1(n+m)\sigma}\phi^{-1}(x)
P^D_{(n+m)\sigma-s}(\phi f)(x)
\end{eqnarray}
and in the last  inequality we used \eqref{uniform boundary2}. It
follows from \eqref{uniform convergence} that there is a constant
$C>0$ such that
\begin{equation}\label{uniform boundary-delta}
p^D(s,y,x)\leq Ce^{\lambda_1 s}\phi(y)\widetilde\phi(x),\qquad
\forall\ s>\sigma,\ x,y\in D.
\end{equation}
Thus
\begin{eqnarray*}
&&\mathbb P_\mu\bigg[\sum_{s\in \widetilde J_{n,m}}
F_{e^{-\lambda_1 (n+m)\sigma}P^D_{(n+m)\sigma-\cdot}(\phi f)}(s, \Delta X_s)^2\bigg]\\
&\leq &C \|f\|_\infty^2\langle
\phi,\mu\rangle\int_D\widetilde\phi(x){\rm d} x\int_{n\sigma}^\infty
e^{-\lambda_1 s} {\rm d}s\int_0^{e^{\lambda_1 s}}r^2 n^\phi(x, {\rm
d}r).
\end{eqnarray*}
Summing over $n$, we get
\begin{eqnarray}\label{I+II}
&&\sum_{n=1}^\infty \mathbb P_\mu\bigg[\sum_{s\in \widetilde J_{n,m}}
F_{e^{-\lambda_1 (n+m)\sigma}P^D_{(n+m)\sigma-\cdot}(\phi f)}(s,
\Delta X_s)^2\bigg]\nonumber\\
&\leq& \sum_{n=1}^\infty  C\| f\|_\infty^2\langle \phi,\mu\rangle
\int_D\widetilde\phi(x){\rm d} x\int_{n\sigma}^\infty e^{-\lambda_1 s}
{\rm d}s\int_0^{e^{\lambda_1 s}}r^2n^\phi(x, {\rm d}r)\nonumber\\
&\leq & C\| f\|_\infty^2\langle
\phi,\mu\rangle\int_D\widetilde\phi(x){\rm d}x\int_0^\infty {\rm d}t
\int_{t\sigma}^\infty  e^{-\lambda_1 s}{\rm d}s\int_0^{
e^{\lambda_1 s}}r^2 n^\phi(x, {\rm d}r)\nonumber\\
&= &\frac{C}{\sigma}\| f\|_\infty^2\langle
\phi,\mu\rangle\int_D\widetilde\phi(x){\rm d}x
\int_{0}^\infty s e^{-\lambda_1 s}{\rm d}s\int_0^{
e^{\lambda_1 s}}r^2 n^\phi(x, {\rm d}r)\nonumber\\
&\le & \frac{C}{\sigma}\| f\|_\infty^2\langle \phi,
\mu\rangle\int_D\widetilde\phi(x){\rm d}x
\int_1^\infty r^2 n^\phi(x, {\rm d}r)\int_{
\lambda_1^{-1}\ln r}^\infty se^{-\lambda_1 s}{\rm d}s\nonumber\\
&+&\frac{C}{\sigma}\| f\|_\infty^2\langle \phi,
\mu\rangle\int_D\widetilde\phi(x){\rm d}x
\int_0^1 r^2 n^\phi(x, {\rm d}r)\int_{ 0}^\infty
se^{-\lambda_1 s}{\rm d}s\nonumber\\
&=:&I+II.
\end{eqnarray}
Using \eqref{assumption on transformed bm} we immediately get that
$II<\infty$. On the other hand,
$$
I=\frac{C}{\lambda_1^2\sigma}\| f\|_\infty^2\langle
\phi,\mu\rangle\int_D\widetilde\phi(x){\rm d}x \int_1^\infty r(\ln
r+1) n^\phi(x, {\rm d}r).
$$
Now we can use $\int_D l(x)\widetilde\phi(x){\rm d}x<\infty$ and
\eqref{assumption on bm} to get that $I<\infty$. The proof of
\eqref{finite 2} is now complete.   For any $\varepsilon>0$,  using
\eqref{finite 2} and Chebyshev's inequality we have
\begin{eqnarray*}
&&\sum_{n=1}^\infty\mathbb P_\mu\left(\big|H_{(n+m)\sigma}(f)- \mathbb
P_\mu[H_{(n+m)\sigma}(f)\big|\mathcal F_{n\sigma}]\big|>\varepsilon\right)\\
&\leq& \varepsilon^{-2}\sum_{n=1}^\infty \mathbb P_\mu\left[H_{(n+m)\sigma}(f) -\mathbb
P_\mu(H_{(n+m)\sigma}(f)\big|\mathcal F_{n\sigma})\right]^2\\
&<& \infty.
\end{eqnarray*}
Then \eqref{limit1} follows easily from the  Borel-Cantelli Lemma.
\qed

\begin{lemma}\label{limit lemma 1(2)}
If $\int_Dl(x)\widetilde\phi(x){\rm d}x<\infty$, then for any $m
\in\mathbb N$, $\sigma>0$, $\mu\in  {\mathcal M}_F(D)^0$ and
$f\in\mathcal B_b^+(D)$ we have
\begin{equation}\label{limit2}
\lim_{n\rightarrow\infty}L_{(n+m)\sigma}(f)-\mathbb P_\mu\left[
L_{(n+m)\sigma}(f)\big|\mathcal F_{n\sigma}\right]=0,\qquad \mathbb
P_\mu{\rm-a.s.}
\end{equation}
\end{lemma}

{\bf Proof:}  It is easy to see that
$$
\mathbb P_\mu\left[ L_{(n+m)\sigma}(f)\big|\mathcal
F_{n\sigma}\right]=e^{-\lambda_1 (n+m)\sigma} \int_0^{n\sigma}\int_D
P^D_{(n+m)\sigma-s}(\phi f)(x) S^{J, (2)}({\rm d}s, {\rm d}x).
$$
Therefore,
\begin{equation}\label{identity 3}
L_{(n+m)\sigma}(f)-\mathbb P_\mu\left[
L_{(n+m)\sigma}(f)\big|\mathcal F_{n\sigma}\right]= e^{-\lambda_1
(n+m)\sigma} \int_{n\sigma}^{(n+m)\sigma}\int_D
P^D_{(n+m)\sigma-s}(\phi f)(x)S^{J, (2)}({\rm d}s, {\rm d}x).
\end{equation}
It follows from Lemma \ref{finite lemma} that, almost surely, the
support of the measure $ N^{(2)}_{\phi} $ consists of finitely many
points. Hence almost surely there exists $N_0\in\mathbb N$ such that
for any $n>N_0$ ,
\begin{eqnarray}\label{equality2}
&&e^{-\lambda_1 (n+m)\sigma} \int_{n\sigma}^{(n+m)\sigma}\int_D
P^D_{(n+m)\sigma-s}(\phi f)(x)
S^{J, (2)}({\rm d}s, {\rm d}x)\nonumber\\
&=& -e^{-\lambda_1 (n+m)\sigma}
\int_{n\sigma}^{(n+m)\sigma}\int_{\mathcal M_F(D)}F_{P^D_{(n+m)
\sigma-\cdot}(\phi f)(\cdot)}(s, \nu)\widehat{N}^{(2)}_\phi({\rm d}
s, {\rm d}\nu)\nonumber\\
&=&-e^{-\lambda_1(n+m)\sigma}\int_{n\sigma}^{(n+m)\sigma}{\rm
d}s\int_D X_s({\rm d}x)\phi(x)^{-1}P^D_{(n+m)\sigma-s}(\phi
f)(x)\int_{e^{\lambda_1 s}}^\infty r n^{\phi}(x, {\rm d}r).
\end{eqnarray}
where in the last equality we used  \eqref{identity1(2)} and
\eqref{formu-F}. Using \eqref{uniform boundary2} we get
\begin{eqnarray}\label{inequality3}
&&\left|e^{-\lambda_1 (n+m)\sigma}
\int_{n\sigma}^{(n+m)\sigma}\int_D P^D_{(n+m)\sigma-s}(\phi f)(x)
S^{J, (2)}({\rm d}s, {\rm d}x)\right|\nonumber\\
&\le&\|f\|_\infty \int_{n\sigma}^\infty e^{-\lambda_1 s}{\rm
d}s\int_D X_s({\rm d}x)\int_{e^{\lambda_1 s}}^\infty  r n^{\phi}(x,
{\rm d}r).
\end{eqnarray}
On the other hand, by \eqref{uniform boundary-delta}, we have
\begin{eqnarray*}
&&\mathbb P_\mu \left[\int_{n\sigma}^\infty e^{-\lambda_1 s}{\rm d}s
\int_D X_s({\rm d}x)\int_{e^{\lambda_1 s}}^\infty  r n^{\phi}(x, {\rm d}r)\right]\\
&=&\int_{n\sigma}^\infty e^{-\lambda_1 s}{\rm d}s \int_D\mu({\rm
d}y)\int_D p^D(s, y, x){\rm d}x\int_{e^{\lambda_1
s}}^\infty rn^\phi(x, {\rm d}r)\\
&\leq& C\langle\phi,\mu\rangle\int_{n\sigma}^\infty{\rm d}s
\int_D\widetilde\phi(x){\rm d}x\int_{e^{\lambda_1 s}}^\infty
r n^{\phi}(x, {\rm d}r)\\
&\leq & C\langle\phi,\mu\rangle\int_D\widetilde\phi(x){\rm d}x
\int_{e^{\lambda_1 n\sigma}}^\infty r n^{\phi}(x, {\rm d}r)
\int_0^{\lambda_1^{-1}\ln r} {\rm d}s\\
&=& \frac{C}{\lambda_1} \langle\phi,
\mu\rangle\int_D\widetilde\phi(x){\rm d}x\int_{e^{\lambda_1
n\sigma}}^\infty  r\ln r n^{\phi}(x, {\rm d}r).
\end{eqnarray*}
Applying the dominated convergence theorem, we obtain that in
$L^1(\mathbb P_\mu)$,
$$
\lim_{n\rightarrow\infty} \int_{n\sigma}^\infty e^{-\lambda_1 s}{\rm
d}s\int_D X_s({\rm d}x)\int_{e^{\lambda_1 s}}^\infty  r n^{\phi}(x,
{\rm d}r)=0.
$$
Since $\int_{n\sigma}^\infty e^{-\lambda_1 s}{\rm d}s\int_D X_s({\rm
d}x)\int_{e^{\lambda_1 s}}^\infty  r n^{\phi}(x, {\rm d}r)$ is
decreasing in $n$, the above limit holds almost surely as well.
Therefore, by \eqref{inequality3}, we have
\begin{equation}\label{limit3}
\lim_{n\rightarrow\infty} e^{-\lambda_1 (n+m)\sigma}
\int_{n\sigma}^{(n+m)\sigma}\int_D (P^D_{(n+m)\sigma-s}\phi
f)(x)S^{J, (2)}({\rm d}s, {\rm d}x)=0,\qquad \mathbb P_\mu{\rm
-a.s.}
\end{equation}
Now \eqref{limit2} follows  from \eqref{identity 3} and
\eqref{limit3}.  The proof is complete.\qed

For any $m,n\in\mathbb N$, $\sigma>0$, set
$$
C_{(n+m)\sigma}(f):=e^{-\lambda_1 (n+m)\sigma}
\int_0^{(n+m)\sigma}\int_D(P^D_{(n+m)\sigma-s}\phi f)(x) S^{C}({\rm
d}s,{\rm d}x), \qquad f\in\mathcal B_b^+(D).
$$
Then $\{C_{n\sigma}(f)\}_{n\in\mathbb N}$ is a martingale with
respect to $(\mathcal F_{n\sigma})$ and
$$
\mathbb P_\mu(C_{(n+m)\sigma}(f)\big|\mathcal
F_{n\sigma})=e^{-\lambda_1 (n+m)\sigma}
\int_0^{n\sigma}\int_D(P^D_{(n+m)\sigma-s}\phi f)(x) S^{C}({\rm d}s,
{\rm d}x).
$$

\begin{lemma}\label{cont-conv}
 For any $m
\in\mathbb N$, $\sigma>0$, $\mu\in  {\mathcal M}_F(D)^0$ and
$f\in\mathcal B_b^+(D)$ we have
\begin{equation}\label{limit5}
\lim_{n\rightarrow\infty}C_{(n+m)\sigma}(f)-\mathbb P_\mu[C_{(n+m)\sigma}(f)
\big|\mathcal F_{n\sigma}]=0,\qquad \mathbb P_\mu {\rm-a.s.}
\end{equation}
\end{lemma}

{\bf Proof:} From the quadratic variation formula \eqref{quadratic variation},
\begin{equation}
\begin{split}
&\mathbb P_\mu\left[H_{(n+m)\sigma}(f)-\mathbb P_\mu(H_{(n+m)\sigma}(f)
\big|\mathcal F_{n\sigma})\right]^2\\
=& \mathbb P_\mu\left[e^{-\lambda_1 (n+m)\sigma}\int_{n\sigma}^{(n+m)
\sigma}\int_D(P^D_{(n+m)\sigma-s}\phi f)(x) S^{C}({\rm d}s,{\rm d}x)\right]^2\\
=&\int_{n\sigma}^{(n+m)\sigma}e^{-2\lambda_1(n+m)\sigma}{\rm d}s\int_D
\mu({\rm d}x)\int_Dp^D(s,x,y)\left(P^D_{(n+m)\sigma-s}\phi f\right)^2(y){\rm d}y \\
=& \int_{n\sigma}^{(n+m)\sigma}e^{-2\lambda_1 s}{\rm d}s\int_D
\mu({\rm d}x)\int_Dp^D(s,x,y)\phi^2(y)\left(P^\phi_{(n+m)\sigma-s}
f\right)^2(y){\rm d}y\\
\le& \|f\|^2_\infty\int_{n\sigma}^{(n+m)\sigma}e^{-\lambda_1 s}{\rm d}s
\int_D\phi(x)\mu({\rm d}x)P^\phi_s\phi(y){\rm d}y\\
\leq & \frac{1}{\lambda_1}\|\phi\|_{\infty}\|f\|^2_\infty\langle
\phi,\mu\rangle e^{-\lambda_1 n\sigma}.
\end{split}
\end{equation}
Therefore, we have
\begin{equation}
\sum_{n=1}^\infty \mathbb P_\mu\left[C_{(n+m)\sigma}(f)-\mathbb
P_\mu(C_{(n+m)\sigma}(f)\big|\mathcal F_{n\sigma})\right]^2<\infty.
\end{equation}
By the Borel-Cantelli lemma, we get \eqref{limit5}.
\qed

Combining the three lemmas above, we  have the following result.

\begin{lemma}\label{limit lemma1}
If $\int_Dl(x)\widetilde\phi(x){\rm d}x<\infty$, then for any $m
\in\mathbb N$, $\sigma>0$, $\mu\in {\mathcal M}_F(D)^0$ and
$f\in\mathcal B_b^+(D)$ we have
\begin{equation}\label{limitmain1}
\lim_{n\rightarrow\infty} e^{-\lambda_1(n+m)\sigma}\langle\phi f,
X_{(n+m)\sigma}\rangle -\mathbb P_\mu\left[
e^{-\lambda_1(n+m)\sigma}\langle\phi f,
X_{(n+m)\sigma}\rangle\big|\mathcal F_{n\sigma}\right]=0,\qquad
\mathbb P_\mu {\rm-a.s.}
\end{equation}
\end{lemma}

{\bf Proof:} From \eqref{dual representation}, we know that
$e^{-\lambda_1(n+m)\sigma}\langle\phi f, X_{(n+m)\sigma}\rangle$ can
be decomposed into three parts:
\begin{eqnarray*}
&&e^{-\lambda_1(n+m)\sigma}\langle\phi f, X_{(n+m)\sigma}\rangle\\
&=&e^{-\lambda_1 (n+m)\sigma}\langle P^D_{(n+m)\sigma}(\phi f),\mu
\rangle+e^{-\lambda_1(n+m)\sigma}
\int_0^{(n+m)\sigma}\int_DP^D_{(n+m)\sigma-s}(\phi f)(x)S({\rm d}s, {\rm d}x)\\
&=&e^{-\lambda_1 (n+m)\sigma}\langle P^D_{(n+m)\sigma}(\phi
f),\mu\rangle+H_{(n+m)\sigma}(f)+L_{(n+m)\sigma}(f)+
 C_{(n+m)\sigma}(f).
\end{eqnarray*}
Therefore,
\begin{eqnarray*}
&&e^{-\lambda_1(n+m)\sigma}\langle\phi f, X_{(n+m)\sigma}\rangle
-\mathbb P_\mu\left[ e^{-\lambda_1(n+m)\sigma}\langle\phi f,
X_{(n+m)\sigma}\rangle\big|\mathcal F_{n\sigma}\right]\\
&=&  H_{(n+m)\sigma}(f)-\mathbb
P_\mu\left[H_{(n+m)\sigma}(f)|\mathcal
F_{n\sigma}\right]+L_{(n+m)\sigma}(f)-\mathbb
P_\mu\left[L_{(n+m)\sigma}(f)|\mathcal F_{n\sigma}\right]\\
  &&+C_{(n+m)\sigma}(f)-\mathbb P_\mu\left[C_{(n+m)\sigma}(f)|\mathcal F_{n\sigma}\right].
\end{eqnarray*}
Now the conclusion of this lemma follows immediately from Lemma
\ref{limit lemma 1}, Lemma \ref{limit lemma 1(2)}, and  Lemma \ref{cont-conv}.\qed

\begin{thm}\label{main skelton}
If $\int_Dl(x)\widetilde\phi(x){\rm d}x<\infty$, then for any
$\sigma>0$, $\mu\in {\mathcal M}_F(D)^0$ and $f\in\mathcal B_b^+(D)$
we have
$$
\lim_{n\rightarrow\infty}e^{-\lambda_1 n\sigma}\langle\phi f,
X_{n\sigma}\rangle=M_{\infty}(\phi)\int_D\widetilde\phi(z)\phi(z)f(z){\rm
d}z ,\qquad \mathbb P_\mu{\rm -a.s.}
$$
\end{thm}

{\bf Proof:} By \eqref{mean-formu} and the Markov property of
super-processes we have
\begin{equation}
\mathbb P_\mu\left[ e^{-\lambda_1(n+m)\sigma}\langle\phi f,
X_{(n+m)\sigma}\rangle\big|\mathcal F_{n\sigma}\right]
=e^{-\lambda_1n\sigma}\langle e^{-\lambda_1m\sigma}
P^D_{m\sigma}(\phi f), X_{n\sigma}\rangle.
\end{equation}
It follows from \eqref{uniform convergence} that
there exist constants $c>0$ and $\nu>0$ such that
$$
\left|\frac{e^{-\lambda_1m\sigma}P^D_{m\sigma}(\phi f)(x)}{\phi(x)}-
\int_D\widetilde\phi(z)\phi(z)f(z){\rm d}z\right|\le c e^{-\nu
m\sigma} \int_D\widetilde\phi(z)\phi(z)f(z){\rm d}z,
$$
which is equivalent to
$$
\left|\frac{e^{-\lambda_1m\sigma}P^D_{m\sigma}(\phi f)(x)}
{\phi(x)\int_D\widetilde\phi(z)\phi(z)f(z){\rm d}z}-1\right|\le c
e^{-\nu m\sigma}.
$$
Thus there exist positive constants $k_m\leq 1$ and $K_m\geq 1$ such
that
$$
k_m\phi(x)\int_D\widetilde\phi(z)\phi(z)f(z){\rm d}z\leq
e^{-\lambda_1m\sigma}P^D_{m\sigma}(\phi f)(x)\leq
K_m\phi(x)\int_D\widetilde\phi(z)\phi(z)f(z){\rm d}z,
$$
and that $\lim_{m\rightarrow\infty}k_m =
\lim_{m\rightarrow\infty}K_m=1$. Hence,
\begin{eqnarray*}
e^{-\lambda_1n\sigma}\langle e^{-\lambda_1m\sigma}P^D_{m\sigma}(\phi
f), X_{n\sigma}\rangle &\geq& k_m e^{-\lambda_1n\sigma}\langle\phi,
X_{n\sigma}\rangle
\int_D\widetilde\phi(z)\phi(z)f(z){\rm d}z\\
&=&k_m M_{n\sigma}(\phi)\int_D\widetilde\phi(z)\phi(z)f(z){\rm
d}z,\qquad \mathbb P_\mu{\rm -a.s.}
\end{eqnarray*}
and
\begin{eqnarray*}
e^{-\lambda_1n\sigma}\langle e^{-\lambda_1m\sigma}P^D_{m\sigma}(\phi
f), X_{n\sigma}\rangle &\leq& K_m e^{-\lambda_1n\sigma}\langle\phi,
X_{n\sigma}\rangle
\int_D\widetilde\phi(z)\phi(z)f(z){\rm d}z\\
&=&K_m M_{n\sigma}(\phi)\int_D\widetilde\phi(z)\phi(z)f(z){\rm d}z,
\qquad \mathbb P_\mu{\rm -a.s.}
\end{eqnarray*}
These two inequalities and Lemma \ref{limit lemma1} imply that
\begin{eqnarray*}
\limsup_{n\rightarrow\infty}e^{-\lambda_1 n\sigma}\langle\phi f,
X_{n\sigma}\rangle &=&
\limsup_{n\rightarrow\infty}e^{-\lambda_1 (n+m)\sigma}\langle
\phi f, X_{(n+m)\sigma}\rangle\\
&= &\limsup_{n\rightarrow\infty}\mathbb P_\mu \left[e^{-\lambda_1
(n+m)\sigma}\langle\phi f, X_{(n+m)\sigma}\rangle\big|\mathcal F_{n\sigma}\right]\\
&=&\limsup_{n\rightarrow\infty}e^{-\lambda_1n\sigma}\langle
e^{-\lambda_1m\sigma}P^D_{m\sigma}(\phi f), X_{n\sigma}\rangle\\
&\leq &\limsup_{n\rightarrow\infty}K_m M_{n\sigma}(\phi)\int_D
\widetilde\phi(z)\phi(z)f(z){\rm d}z\\
&=& K_m M_{\infty}(\phi)\int_D\widetilde\phi(z)\phi(z)f(z){\rm
d}z,\qquad \mathbb P_\mu{\rm -a.s.}
\end{eqnarray*}
and that
$$
\liminf_{n\rightarrow\infty}e^{-\lambda_1 n\sigma}\langle\phi f,
X_{n\sigma}\rangle \geq k_m M_{\infty}(\phi)\int_D
\widetilde\phi(z)\phi(z)f(z){\rm d}z,\qquad \mathbb P_\mu{\rm -a.s.}
$$
Letting $m\to\infty$, we get
$$
\lim_{n\rightarrow\infty}e^{-\lambda_1 n\sigma}\langle\phi f,
X_{n\sigma}\rangle = M_{\infty}(\phi)
\int_D\widetilde\phi(z)\phi(z)f(z){\rm d}z,\qquad \mathbb P_\mu{\rm
-a.s.}
$$
The proof is now complete.\qed

We are now ready to give the proof of Theorem \ref{main theorem}.

\vspace{.1in}

{\bf Proof of Theorem  \ref{main theorem}}\quad Put
$\Delta_{\sigma}(f):=\sup_{0\leq t\leq \sigma}
\|P^\phi_{t}f-f\|_{\infty}$. Then for $t\in[n\sigma, (n+1)\sigma]$,
\begin{eqnarray}\label{inequality 4}
\left|e^{-\lambda_1 t}\left\langle\phi P_{(n+1)\sigma-t}^\phi
f,X_t\right\rangle- e^{-\lambda_1 t}\langle\phi f,X_t\rangle\right| \leq
e^{-\lambda_1 t}\left\langle\phi|P_{(n+1)\sigma-t}^\phi f-f|,X_t\right\rangle
\leq M_t(\phi)\Delta_{\sigma}(f).
\end{eqnarray}
By the strong continuity of the semigroup $(P_t^\phi)$ in
$L^\infty(D)$, we have $\lim_{\sigma\rightarrow
0}\Delta_{\sigma}(f)=0$. Thus,
\begin{eqnarray}\label{limit-0}
\lim_{\sigma\to 0}\lim_{n\to\infty}\sup_{t\in[n\sigma,
(n+1)\sigma]}\left|e^{-\lambda_1 t}\left\langle\phi
P_{(n+1)\sigma-t}^\phi f,X_t\right\rangle- e^{-\lambda_1
t}\left\langle\phi f,X_t\right\rangle\right|=0,\quad\mathbb
P_\mu{\rm -a.s.}
\end{eqnarray}
Therefore, to prove Theorem \ref{main theorem}, we only need to show
that
\begin{eqnarray}\label{toprove-limit}
\lim_{\sigma\to 0}\lim_{n\to\infty}\sup_{t\in[ n\sigma,\
(n+1)\sigma]}e^{-\lambda_1 t}\langle\phi P_{(n+1)\sigma-t}^\phi
f,X_t\rangle=M_\infty(\phi)\int_D\phi(x)\widetilde\phi(x) f(x){\rm
d}x,\quad \mathbb P_\mu{-a.s.}
\end{eqnarray}
For any $n\in\mathbb N$ and $\sigma>0$, $(X_t, t\in[n\sigma,
(n+1)\sigma], \mathbb P_\mu(\cdot|\mathcal F_{n\sigma}))$ can be
regarded as a $(\xi^D, \psi(\lambda)-\beta\lambda)$-super-diffusion
with initial value $X_{n\sigma}$.  Thus, for arbitrary $g\in\mathcal
B_b^+(D)$, we have by \eqref{dual representation}
\begin{eqnarray*}
e^{-\lambda_1 t}\langle\phi g,X_t\rangle = e^{-\lambda_1 t}\langle
P^D_{t-n\sigma}(\phi g),X_{n\sigma}\rangle+e^{-\lambda_1
t}\int_{n\sigma}^{t}\int  P^D_{t-s}(\phi g)(x)S({\rm d}s, {\rm
d}x),\qquad t\in[n\sigma, (n+1)\sigma].
\end{eqnarray*}
Taking $g(x)=P_{(n+1)\sigma-t}^\phi f(x)$ in the above identity  and
using \eqref{def-p-phi}, we get
\begin{eqnarray}\label{truncate type identity}
e^{-\lambda_1 t}\left\langle\phi P_{(n+1)\sigma-t}^\phi f,X_t\right\rangle =
e^{-\lambda_1 n\sigma}\left\langle \phi P^\phi_\sigma
f,X_{n\sigma}\right\rangle+\int_{n\sigma}^{t}e^{-\lambda_1 s}\int ( \phi
P^\phi_{(n+1)\sigma-s}f)(x)S({\rm d}s,{\rm d}x).
\end{eqnarray}
Since $\widetilde\phi\phi$ is the invariant probability density of
the semigroup $(P^\phi_t)$, we have by Lemma \ref{main skelton},
\begin{eqnarray}\label{first-term-limit}
\lim_{n\to\infty}e^{-\lambda_1 n\sigma}\langle \phi P^\phi_\sigma f,
X_{n\sigma}\rangle&=&M_\infty(\phi)\int_D\phi(x)\widetilde\phi(x)
P^\phi_\sigma f(x){\rm d}x\nonumber\\
&=&M_\infty(\phi)\int_D\phi(x)\widetilde\phi(x) f(x){\rm d}x.
\end{eqnarray}
Hence, by \eqref{truncate type identity} and
\eqref{first-term-limit}, to prove \eqref{toprove-limit} it suffices
to show that
\begin{eqnarray}\label{toprove-limit2}
\lim_{\sigma\to 0}\lim_{n\to\infty}\sup_{t\in[n\sigma, (n+1)\sigma]}
\int_{n\sigma}^{t}e^{-\lambda_1 s}\int ( \phi
P^\phi_{(n+1)\sigma-s}f)(x)S({\rm d}s,{\rm d}x)=0,\quad \mathbb
P_\mu{\rm -a.s.}
\end{eqnarray}
Since $S({\rm d}s,{\rm d}x)=S^J({\rm d}s, {\rm d}x)+S^C({\rm d}s,
{\rm d}x)=S^{J, (1)}({\rm d}s, {\rm d}x)+S^{J, (2)}({\rm d}s, {\rm
d}x)+S^C({\rm d}s, {\rm d}x),$ we have
\begin{eqnarray*}
&&\int_{n\sigma}^{t}e^{-\lambda_1 s}\int ( \phi
P^\phi_{(n+1)\sigma-s}f)(x)S({\rm d}s,{\rm d}x)\\
&=&\int_{n\sigma}^{t}e^{-\lambda_1 s}\int ( \phi
P^\phi_{(n+1)\sigma-s}f)(x)S^{J, (1)}({\rm d}s, {\rm d}x)+
\int_{n\sigma}^{t}e^{-\lambda_1 s}\int ( \phi
P^\phi_{(n+1)\sigma-s}f)(x)S^{J, (2)}({\rm d} s, {\rm d}x)\\
&&+\int_{n\sigma}^{t}e^{-\lambda_1 s}\int ( \phi
P^\phi_{(n+1)\sigma-s}f)(x)S^{C}({\rm d} s, {\rm d}x)\\
&=:&H_{n,t}^\sigma(f)+L_{n,t}^\sigma( f)+C_{n,t}^\sigma(f).
\end{eqnarray*}
Thus we only need to prove that
\begin{equation}\label{H-zero-limit}
\lim_{\sigma\to 0}\lim_{n\to\infty}\sup_{t\in[n\sigma,
(n+1)\sigma]}H_{n,t}^\sigma(f)=0, \qquad \mathbb P_\mu{\rm -a.s.},
\end{equation}
\begin{equation}\label{L-zero-limit}
\lim_{\sigma\to 0}\lim_{n\to\infty}\sup_{t\in[n\sigma,
(n+1)\sigma]}L_{n,t}^\sigma( f)=0,\qquad \mathbb P_\mu{\rm -a.s.},
\end{equation}
and
\begin{equation}\label{C-zero-limit}
\lim_{\sigma\to 0}\lim_{n\to\infty}\sup_{t\in[n\sigma,
(n+1)\sigma]}C_{n,t}^\sigma( f)=0,\qquad \mathbb P_\mu{\rm -a.s.}
\end{equation}

 It follows from Chebyshev's inequality that, for any
$\varepsilon>0$, we have
\begin{eqnarray}\label{inequality 3}
\mathbb P_\mu\left(\sup_{t\in[n\sigma,
(n+1)\sigma]}\left|H_{n,t}^\sigma(f)\right|>\varepsilon\right)\leq
\frac{1}{\varepsilon^2}\mathbb P_\mu\left(\sup_{t\in[n\sigma,
(n+1)\sigma]}\int_{n\sigma}^t e^{-\lambda_1 s}\int\left( \phi
P^\phi_{(n+1)\sigma-s}f\right)(x)S^{J, (1)}({\rm d} s, {\rm d}x)
\right)^2.
\end{eqnarray}
Since the process $(H_{n,t}^\sigma(f); t\in[n\sigma, (n+1)\sigma])$
is a martingale with respect to $(\mathcal F_t)_{t\in[n\sigma,
(n+1)\sigma]}$, applying Burkholder-Davis-Gundy inequality to
$H_{n,t}^\sigma(f)$ and using an argument similar to the one in the
proof of Lemma \ref{limit lemma 1}, we obtain
\begin{eqnarray}\label{I-1-2}
&&\mathbb P_\mu\Big(\sup_{t\in[n\sigma, (n+1)\sigma]}\int_{n\sigma}^t
e^{-\lambda_1 s}\int ( \phi P^\phi_{(n+1)\sigma-s}f)(x)S^{J, (1)}
({\rm d} s, {\rm d}x)\Big)^2\nonumber\\
&\leq &C_1\mathbb P_\mu\Big(\int_{n\sigma}^{(n+1)\sigma}
e^{-\lambda_1 s}\int\phi(x) P^\phi_{(n+1)\sigma-s}f(x)
S^{J, (1)}({\rm d} s, {\rm d}x)\Big)^2\nonumber\\
&=&C_1\int_{n\sigma}^{(n+1)\sigma}e^{-2\lambda_1 s}{\rm d}s
\int_D\mu({\rm d} y)\int_D{\rm d}xp^D(s, y, x)
\left(P^\phi_{(n+1)\sigma-s}f\right)^2(x)\int_0^{
e^{\lambda_1 s}}r^2n^\phi(x, {\rm d}r)\nonumber\\
&=&C_1\int_{n\sigma}^{(n+1)\sigma}e^{-\lambda_1 s}{\rm d}s
\int_D\phi(y)\mu({\rm d} y)\int_D{\rm d}xp^\phi(s, y, x)
\phi(x)^{-1}\left(P^\phi_{(n+1)\sigma-s}f\right)^2(x)\int_0^{
e^{\lambda_1 s}}r^2n^\phi(x, {\rm d}r)\nonumber\\
&\leq&C_2(\sigma)\|f\|^2_{\infty}\langle\phi,\mu\rangle
\int_{n\sigma}^{(n+1)\sigma}e^{-\lambda_1 s}{\rm d}s \int_D{\rm
d}x\widetilde\phi(x)\int_0^{e^{\lambda_1 s}}
r^2n^\phi(x, {\rm d}r)\nonumber\\
&\leq&C_2(\sigma) \|f\|_{\infty}^2\langle\phi,\mu\rangle
\left[\int_{n\sigma}^{(n+1)\sigma}e^{-\lambda_1 s}{\rm d}s
\int_D\widetilde\phi(x){\rm d}x\left\|\int_0^1r^2n^\phi(x,
{\rm d}r)\right\|_\infty\right.\nonumber\\
&&+\left.\int_{n\sigma}^{(n+1)\sigma}e^{-\lambda_1 s}{\rm d}s
\int_D\widetilde\phi(x){\rm d}x\int_1^{e^{\lambda_1 s}}r^2
n^\phi(x, {\rm d}r)\right]\nonumber\\
&=:& C_2(\sigma)
\|f\|_{\infty}^2\langle\phi,\mu\rangle(I_1(n)+I_2(n)),
\end{eqnarray}
where  $C_1$ and  $C_2(\sigma)$ are positive constants independent
of $n$. \eqref{assumption on bm} implies that
$\|\int_0^1r^2n^\phi(\cdot, {\rm d}r)\|_\infty<\infty$. Thus
\begin{equation}\label{I-1}
\sum^\infty_{n=1}I_1(n)<\infty.
\end{equation}
Using Fubini's theorem, we have
\begin{eqnarray}\label{I-2}
\sum_{n=1}^\infty I_2(n)&=&\sum_{n=1}^\infty\int_{n\sigma}^{(n+1)\sigma}
e^{-\lambda_1 s}{\rm d}s\int_D\widetilde\phi(x){\rm d}x\int_1^{
e^{\lambda_1 s}}r^2n^\phi(x, {\rm d}r)\nonumber\\
&=&\int_\sigma^\infty e^{-\lambda_1 s}{\rm d}s\int_D
\widetilde\phi(x){\rm d}x\int_1^{e^{\lambda_1 s}}r^2
n^\phi(x, {\rm d}r)\nonumber\\
&\leq&\int_D\widetilde\phi(x){\rm d}x\int_1^\infty r^2
n^\phi(x, {\rm d}r)\int_{\lambda_1^{-1}\ln r}^\infty
e^{-\lambda_1 s}{\rm d}s\nonumber\\
&=&\lambda_1^{-1}\int_D\widetilde\phi(x){\rm d}x
\int_1^\infty rn^\phi(x, {\rm d}r)\nonumber\\
&\le &\lambda_1^{-1}\int_D\widetilde\phi(x){\rm d}x
\left\|\int_1^\infty r n^\phi(x, {\rm d}r)\right\|_\infty <\infty.
\end{eqnarray}
Combining \eqref{inequality 3}, \eqref{I-1-2}, \eqref{I-1} and
\eqref{I-2}, we get that, for any $\varepsilon>0$,
\begin{equation}
\sum_{n=1}^\infty\mathbb P_\mu\left(\sup_{t\in[n\sigma,
(n+1)\sigma]}\left|H_{n,t}^\sigma(f)\right|>\varepsilon\right)<\infty.
\end{equation}
Thus by the Borel-Cantelli lemma we have, for any $\sigma>0$,
\begin{equation}\label{limit4}
\lim_{n\rightarrow\infty}\sup_{t\in[n\sigma, (n+1)\sigma]}
H_{n,t}^\sigma(f)=0,\qquad  \mathbb P_\mu{\rm -a.s.}
\end{equation}
Therefore \eqref{H-zero-limit} is valid.

Similarly, we can prove that
\begin{equation}\label{limit6}
\lim_{n\rightarrow\infty}\sup_{t\in[n\sigma, (n+1)\sigma]}
|C_{n,t}^\sigma(f)|=0,\qquad  \mathbb P_\mu{\rm -a.s.},
\end{equation}
and then we get \eqref{C-zero-limit}. The details are omitted here.

Using an argument similar to \eqref{equality2}
we can see that almost surely there exists $N_0\in\mathbb N$ such
that when $n>N_0$,
\begin{eqnarray}\label{inequality 5}
\left|L_{n,t}^\sigma(f)\right|&=&\int_{n\sigma}^te^{-\lambda_1s}
{\rm d}s\int_D P^\phi_{(n+1)\sigma-s}f(x)X_s({\rm d}x)
\int_{e^{\lambda_1 s}}^\infty r n^\phi(x,{\rm d}r)\nonumber\\
&\leq& \|f\|_{\infty}\int_{n\sigma}^{(n+1)\sigma}
e^{-\lambda_1 s}X_s({\rm d}x)\phi(x)\int_{e^{\lambda_1 s}}^\infty
rn^\phi(x, {\rm d}r)\nonumber\\
&\leq& \sigma\|f\|_{\infty}\left\|\int_1^\infty r
n^\phi(x, {\rm d}r)\right\|_{\infty}\sup_{s\in[n\sigma, (n+1)
\sigma]}M_s(\phi).
\end{eqnarray}
Therefore \eqref{L-zero-limit} holds. The proof of Theorem \ref{main
theorem} is now complete.

The following result strengthens Theorem \ref{main theorem} in the
sense that the exceptional does not depend on $f$ and $\mu$.

\begin{thm}\label{main theorem2}
Suppose that Assumption \ref{assume1} holds, $\lambda_1>0$ and that
$X$ is a $(\xi^D, \psi(\lambda)-\beta\lambda)-$super-diffusion. Then there
exists $\Omega_0\subset\Omega$ of probability one (that is,
$\mathbb P_\mu (\Omega_0)=1$ for every $\mu\in \mathcal M_F(D)^0$) such that, for
every $\omega\in\Omega_0$ and for every bounded Borel measurable
function $f$ on $\R^d$ with compact support whose set of
discontinuous points has zero Lebesgue measure, we have
\begin{eqnarray}\label{main result2}
\lim_{t\rightarrow\infty}e^{-\lambda_1 t}\langle f, X_t\rangle =
M_\infty(\phi)\int_D\widetilde\phi(y)f(y){\rm d}y.
\end{eqnarray}
\end{thm}

{\bf Proof:}  Note that there exists a countable base ${\cal U}$ of
open sets $\{U_k, k\ge 1\}$ that is closed under finite unions.
Define
$$
\Omega_0:= \left\{\omega \in \Omega: \,
\lim_{t\rightarrow\infty}e^{-\lambda_1 t}\langle I_{U_k} \phi,
X_t\rangle = M_\infty(\phi)\int_{U_k}\phi(y)\widetilde\phi(y){\rm
d}y\quad \mbox{for every }k\ge 1 \right\}.
$$
By Theorem \ref{main theorem}, for any $\mu\in \mathcal M_F(D)^0$,
$P_\mu(\Omega_0)=1$.  For any open set $U$, there exists a sequence
of increasing open sets $\{U_{n_k}; k\ge 1\}$ in ${\cal U}$ so that
$\bigcup^\infty_kU_{n_k}=U$. Then for every $\omega\in\Omega_0$,
$$
\liminf_{t\rightarrow\infty}e^{-\lambda_1 t}\langle I_{U} \phi,
X_t\rangle \ge \lim_{t\rightarrow\infty}e^{-\lambda_1 t}\langle
I_{U_{n_k}}\phi, X_t\rangle\ge
M_\infty(\phi)\int_{U_{n_k}}\phi(y)\widetilde\phi(y){\rm
d}y\quad\mbox{ for every }k\ge 1.
$$
Letting $k\to\infty$ yields
\begin{eqnarray}\label{limit-openset}
\liminf_{t\rightarrow\infty}e^{-\lambda_1 t}\langle I_U\phi,
X_t\rangle \ge M_\infty(\phi)\int_U\phi(y)\widetilde\phi(y){\rm
d}y,\quad \mathbb P_\mu{\rm -a.s.} \mbox{ for any }\mu \in \mathcal
M_F(D)^0.
\end{eqnarray}

We now consider \eqref{main result2} on $\{M_\infty(\phi)>0\}$.
For each $\omega\in\Omega_0\cap \{M_\infty(\phi)>0\}$ and $t\ge 0$,
we define two probability measures $\nu_t$ and $\nu$ on $D$
respectively by
$$
\nu_t(A)(\omega)=\frac{e^{\lambda_1t}\langle I_A\phi,
X_t\rangle(\omega)}{M_\infty(\phi)(\omega)},\quad \mbox{and}\quad
\nu(A)=\int_A\phi(y)\widetilde \phi(y)dy,\quad A\in{\cal B}(D).
$$
Note that the measure $\nu_t$ is well-defined for every $t\ge 0$.
 \eqref{limit-openset} tells us that $\nu_t$ converges weakly to
$\nu$ as $t\to\infty$. Since $\phi$ is strictly positive and
continuous on $D$, for every function $f$ on $D$ with compact
support on E whose discontinuity set has zero Lebesgue-measure
(equivalently zero $\nu$-measure), $g:= f/\phi$ is a bounded
function with compact support with the same set of discontinuity. We
thus have
$$
\int_Dg(x)\nu_t(dx)=\int_Dg(x)\nu(dx),
$$
which is equivalent to say
\begin{eqnarray}\label{limit-positive M}
\lim_{t\rightarrow\infty}e^{-\lambda_1 t}\langle f, X_t\rangle =
M_\infty(\phi)\int_D\widetilde\phi(y)f(y){\rm d}y,\qquad
\mbox{for every }\omega\in\Omega_0\cap \{M_\infty(\phi)>0\}.
\end{eqnarray}

Since
$$
e^{-\lambda_1 t}\left|\langle f, X_t\rangle \right| \le
e^{-\lambda_1 t}\langle \left|f\right| , X_t\rangle = e^{-\lambda_1
t}\langle \left|g\phi\right| , X_t\rangle\le \|g\|_\infty
M_\infty(\phi).
$$
\eqref{limit-positive M} holds automatically on
$\{M_\infty(\phi)=0\}$. This completes the proof of the theorem.

\begin{lemma}\label{mean-conv}
For any $f\in {\cal B}_b(D)$, $\mu\in \mathcal M_F(D)^0$,
\begin{equation}\label{mean-limit}
\lim_{t\to\infty} e^{-\lambda_1t}\mathbb P_\mu\langle f,
X_t\rangle=\langle\phi, \mu\rangle\int_Df(y)\widetilde\phi(y)dy.
\end{equation}
\end{lemma}

{\bf Proof:} It follows from \eqref{mean-formu} that
\begin{eqnarray*}
e^{-\lambda_1t}\mathbb P_\mu\langle f,
X_t\rangle&=&\int_D\mu(dx)e^{-\lambda_1t}P^D_tf(x)\\
&=&
\int_D\mu(dx)\int_De^{-\lambda_1t}p^D(t,x,y)f(y)dy\\
&=&\int_D\mu(dx)\phi(x)\int_Dp^\phi(t,x,y)\frac{f(y)}{\phi(y)}dy.
\end{eqnarray*}
Using \eqref{uniform convergence} and the dominated convergence
theorem, we get \eqref{mean-limit}.

\vspace{.1in}

{\bf Proof of Theorem  \ref{new main theorem}}\quad Theorem \ref{new
main theorem} is simply a combination of Theorem \ref{main theorem2}
and Lemma \ref{mean-conv}.

\qed

\begin{singlespace} 
\end{singlespace}
\end{doublespace}

\vskip 0.3truein \vskip 0.3truein

 \noindent {\bf Rong-Li Liu:} Department of Mathematics, Nanjing
 University,  Nanjing, 210093, P. R. China,
  E-mail: {\tt rlliu@nju.edu.cn} \\

\bigskip
 \noindent {\bf Yan-Xia Ren:} LMAM School of Mathematical Sciences,  Center for Statistical Science, Peking
 University,  Beijing, 100871, P. R. China,
  E-mail: {\tt yxren@math.pku.edu.cn} \\

\bigskip
\noindent {\bf Renming Song:} Department of Mathematics,
 University of Illinois,  Urbana, IL 61801 U.S.A.,
  E-mail: {\tt rsong@math.uiuc.edu} \\
 \bigskip
\end{document}